# A Restricted Chen-Nagano Variational Principle for the Einstein-Hilbert Functional

Sergey Stepanov, Irina Tsyganok

**Abstract.** This paper introduces a restricted Chen-Nagano variational principle for the well-known Einstein-Hilbert functional

$$S(g) = \int_M s_g \, dv_g,$$

where $s_g$ is the scalar curvature of the metric $g$ of a compact $n$-dimensional Riemannian manifold $(M, g)$. Instead of allowing arbitrary symmetric variations of the metric, we restrict the admissible variational directions to the infinite-dimensional space of Chen-Nagano harmonic tensors. We derive the corresponding restricted Euler-Lagrange equation and prove that $g$ is critical with respect to all volume-preserving Chen-Nagano variations if and only if there exists a smooth 1-form $\theta$ such that

$$E_g = B_g^*(\theta) - \frac{n-2}{2n} \bar{s}_g \, g,$$

where $B_g^*$ denotes the first-order differential operator, $E_g$ is the Einstein tensor and $\bar{s}_g$ is the average scalar curvature of $g$. This equation may be regarded as a restricted counterpart of the classical Einstein equation. We further show that Chen-Nagano critical metrics naturally give rise to generalized Ricci almost soliton-type geometries and, in the gradient case, to gradient Ricci almost solitons. Moreover, the gauge 1-form $\theta$ is automatically in the kernel of the Sampson Laplacian; equivalently, $\theta^\sharp$ is an infinitesimal harmonic transformation. Finally, we establish global rigidity results. In particular, every compact Chen-Nagano critical metric of constant scalar curvature is Einstein. Further rigidity follows under quasi-negative Ricci curvature and, in the gradient case, under a natural scalar-curvature monotonicity condition along the potential field. Finally, we study the restricted second variation at negative Einstein metrics.

## 1. Introduction

Let $(M, g)$ be a compact Riemannian manifold of dimension $n \geq 3$. The *Einstein-Hilbert functional*

$$S(g) = \int_M s_g \, dv_g,$$

where $s_g$ is the scalar curvature of the metric $g$, occupies a central position in Riemannian geometry and mathematical relativity.

One of the fundamental results of the variational theory states that Einstein metrics arise as critical points of this functional on the space of Riemannian metrics of fixed volume. More precisely, if admissible variations are allowed to range over the entire Fréchet space $C^\infty(S^2M)$ of smooth symmetric bilinear forms, then the corresponding Euler-Lagrange equation is precisely the Einstein equation (see, e.g., Besse [1, Chapter 4, §4.17]). Background material on the variational theory of the Einstein-Hilbert functional can be found in the classical works of Berger and Ebin [2], Koiso [8], Schoen [14], and in the monograph of Besse [1].

The present paper is motivated by the following natural question:

To what extent can the class of admissible metric variations be restricted while preserving a meaningful geometric variational theory?

In this work we investigate this question within the framework of the Chen-Nagano gauge condition

$$\delta\, h = -\frac{1}{2} d\left(\mathrm{tr}_g h\right)$$

where $\delta: C^\infty(S^2M) \to C^\infty(T^*M)$ is the divergence operator for $h \in C^\infty(S^2M)$, introduced by Chen and Nagano in their theory of harmonic metrics and tensors (see [3] and [4]). Symmetric tensor fields satisfying this condition will be called *Chen-Nagano harmonic tensors*, and the corresponding metric deformations will be referred to as Chen-Nagano variations.

The Chen-Nagano condition and related Berger-Ebin techniques have appeared in several geometric contexts, including harmonic maps (in particular, harmonic diffeomorphisms and harmonic submersions), infinitesimal harmonic transformations, and Ricci almost solitons (see [11]; [12]; [15] and [16]).

In the context of geometric evolution equations, the behavior of scalar and Ricci curvatures has been extensively studied from both dynamic and stationary perspectives. In particular, Rovenski, Stepanov, and Tsyganok [13] investigated rigidity phenomena for Hamilton's Ricci flow and Ricci solitons, establishing conditions under which these structures reduce to Einstein or Ricci-flat geometries. While those results were obtained in the setting of curvature evolution equations, the present paper approaches similar rigidity questions from a variational viewpoint by restricting the admissible metric deformations through the Chen-Nagano gauge condition.

It is worth noting that the differential operator $B_g := \delta\, h + \frac{1}{2} d\big(\mathrm{tr}_g h\big)$ considered in this paper coincides with the operator that later became known in Ricci flow theory (see, for example, [6]) as the *Bianchi operator* (see [9]). However, our approach is motivated by the variational theory of *Chen-Nagano harmonic tensors* rather than by gauge fixing in geometric evolution equations.

Relatedly, recent developments in the theory of critical point metrics continue to illustrate the effectiveness of constrained variational methods in global Riemannian geometry under some Ricci curvature constraints (see, e.g., Li and Yu [19]).

The motivation for such a restriction comes from a familiar phenomenon in variational calculus. In many geometric and analytical problems, constrained variational principles often reveal structures that remain hidden in the unrestricted setting. The central question is therefore whether a geometrically natural constraint can still retain enough information to produce meaningful Euler–Lagrange equations. The present paper shows that the *Chen-Nagano gauge* provides such a constraint in the context of the Einstein-Hilbert functional.

The core observation of the paper is that restricting the admissible variations to the infinite-dimensional space of Chen-Nagano tensors does not destroy the geometric content of the Einstein-Hilbert variational principle. On the contrary, the resulting constrained variational problem still leads to a strong Euler-Lagrange equation and yields significant geometric consequences.

From an analytical point of view, the construction may be viewed as an infinite-dimensional analogue of a constrained variational problem. The role of the constraint is played by the Chen-Nagano gauge condition, while the corresponding admissible variations form a proper infinite-dimensional subspace of the space of all symmetric tensor fields (see [3] and [4]).

To formulate the main result, let $B_g^*$ denote the first-order overdetermined elliptic operator

$$B_g^*(\theta) = \delta^*\theta + \frac{1}{2}(\delta\theta)g,$$

where $\delta^*\theta = \frac{1}{2}\mathcal{L}_X g$ for the Lie derivative $\mathcal{L}_X g$ of the metric tensor $g$ along the vector field $X = \theta^\#$ which is dual (by $g$) to the 1-form $\theta$. Its $L^2$-formal adjoint operator is given by

$$B_g(h) = \delta h + \frac{1}{2} d\big(\mathrm{tr}_g h\big).$$

The kernel $\mathcal{H}_g = \mathrm{Ker}\, B_g$ coincides with the infinite-dimensional space of harmonic Chen-Nagano tensors. We can now state the main result, which gives a restricted analogue of Einstein's classical variational principle. Let $E_g = Ric_g - \frac{1}{2}\, s_g\, g$ be the standard *Einstein tensor*. Then the following theorem holds.

**Theorem 3.1.** *Let $(M, g)$ be a compact Riemannian manifold of dimension $n \geq 3$. Then $g$ is a critical point of the Einstein-Hilbert functional with respect to all volume-preserving Chen-Nagano variations if and only if there exists a smooth 1-form $\theta$ such that:*

$$E_g = B_g^*(\theta) - \frac{n-2}{2n}\bar{s}\, g.$$

*where $\bar{s}$ is the average scalar curvature of $(M, g)$.*

The restricted Chen-Nagano variational principle considered here should be viewed in the broader context of constrained variational problems in Riemannian geometry. Classical examples include the theory of critical point metrics, where the Einstein-Hilbert functional is restricted to the space of metrics with constant scalar curvature (see Besse [1] and the survey article of Schoen [14]). In the present work, the constraint is of a different nature and is imposed on the admissible variations rather than on the ambient space of metrics. This leads to a distinct Euler–Lagrange equation and to a different geometric framework.

Thus, the classical Einstein equation is replaced by a natural structural equation involving the Chen-Nagano gauge data. We therefore obtain a restricted Chen-Nagano variational principle for the Einstein-Hilbert functional.

The geometric consequences of equation (3.4) are developed in the subsequent two sections. In Section 4, we relate Chen-Nagano critical metrics to generalized Ricci almost soliton-type geometries and, in the gradient case, to gradient Ricci almost solitons. In Section 5, we prove that the gauge 1-form associated with a Chen-Nagano critical metric is automatically an infinitesimal harmonic transformation and establish several global rigidity results. In particular, every compact Chen-Nagano critical metric of constant scalar curvature is Einstein. Further rigidity is obtained under quasi-negative Ricci curvature and, in the gradient case, under a natural monotonicity condition on the scalar curvature along the potential field.

The paper is organized as follows. Section 2 contains the necessary background on the operators $B_g$ and $B_g^*$, Chen-Nagano harmonic tensors, and the relevant orthogonal decompositions of the space of symmetric 2-tensors. It also establishes a relation between the Chen-Nagano condition, transverse-traceless tensors, and infinitesimal harmonic transformations. Section 3 develops the restricted variational principle and proves the main theorem. Section 4 is devoted to the geometric interpretation of Chen-Nagano critical metrics and their relation to generalized Ricci almost solitons. Finally, Section 5 proves that the associated gauge 1-form is an infinitesimal harmonic transformation and establishes global rigidity results for compact Chen-Nagano critical metrics.
Finally, Section 6 studies the restricted second variation at negative Einstein metrics. We prove that every infinitesimally scalar-curvature-preserving Chen-Nagano variation is transverse-traceless. Consequently, the restricted stability problem reduces exactly to the classical spectral problem for the Einstein operator on $TT$-tensors.

## 2. Chen-Nagano harmonic tensors and their $L^2$-orthogonal decomposition

Let $(M, g)$ be a compact Riemannian manifold of dimension $n \geq 3$. Denote by $C^\infty(S^2M)$ the Fréchet space of smooth symmetric covariant 2-tensor fields on $M$, and by $C^\infty(T^*M)$ the space of smooth 1-forms.

For tensor $h, k \in C^\infty(S^2T^*M)$, their local components in a coordinate system $x^1, \dots, x^n$ will be denoted by $h = h_{ij}\, dx^i \otimes dx^j$ and $k = k_{ij}\, dx^i \otimes dx^j$ The metric $g$ induces the natural $L^2$-inner product

$$\langle h, k\rangle_{L^2} = \int_M g^{ik}\, g^{jl} h_{ij} k_{kl}\, dv_g,$$

where $(g^{ij})$ are the contravariant components of the inverse metric tensor $g^{-1}$.

Let $\delta: C^\infty(S^2M) \to C^\infty(T^*M)$ be the divergence operator locally given by

$$(\delta h)_i = -\nabla^j h_{ji},$$

where $\nabla^i = g^{ij}\nabla_j$, and $\nabla_j$ denotes the covariant derivative with respect to the coordinate vector field $\frac{\partial}{\partial x^j}$. Then the classical Killing operator $\delta^*: C^\infty(T^*M) \to C^\infty(S^2M)$ be its $L^2$-formal adjoint operator defined by (see [2])

$$(\delta^*\theta)_{ij} = \frac{1}{2}(\nabla_i\theta_j + \nabla_j\theta_i).$$

where $\theta_j = \theta\left(\frac{\partial}{\partial x^j}\right)$ are the local covariant components of the 1-form $\theta$.

Define the first-order differential operator $B_g^*: C^\infty(T^*M) \to C^\infty(S^2M)$ by

$$B_g^*(\theta) = \delta^*\theta + \frac{1}{2}(\delta\theta)g.$$

To determine the formal $L^2$-adjoint

$$B_g: C^\infty(S^2M) \to C^\infty(T^*M),$$

consider arbitrary tensors

$$h \in C^\infty(S^2M), \qquad \theta \in C^\infty(T^*M).$$

Using integration by parts, we obtain

$$\langle B_g^*(\theta), h\rangle_{L^2} = \langle \delta^*\theta, h\rangle_{L^2} + \frac{1}{2}\int_M (\delta\theta)\,(\mathrm{tr}_g\, h)\, dv_g,$$

where $\delta\theta = -\nabla^i\theta_i$ is the divergence of the 1-form $\theta$.

Since

$$\langle \delta^*\theta, h\rangle_{L^2} = \langle \theta, \delta h\rangle_{L^2},$$

and

$$\langle \delta\theta, f\rangle_{L^2} = \langle \theta, df\rangle_{L^2}$$

for every smooth function $f \in C^\infty(M)$, we obtain

$$\frac{1}{2}\int_M (\delta\theta)(\operatorname{tr}_g h)\, dv_g = \frac{1}{2}\langle \theta, d(\operatorname{tr}_g h)\rangle_{L^2}.$$

Therefore,

$$\langle B_g^*(\theta), h\rangle_{L^2} = \langle \theta, \delta h + \frac{1}{2} d(\operatorname{tr}_g h)\rangle_{L^2}.$$

Hence, the formal $L^2$-adjoint of $B_g^*$ is given by

$$B_g(h) = \delta h + \frac{1}{2} d(\operatorname{tr}_g h).$$

**Remark 2.1.** The operator $B_g$ is closely related to the classical Bianchi operator appearing in Ricci flow theory and in gauge-fixed formulations of the Einstein equations. In particular, the operator $B_g(h) = \delta h + \frac{1}{2} d(\operatorname{tr}_g h)$ coincides with the coordinate-free form of the Bianchi operator used in modern geometric analysis.

Historically, however, the geometric origin of this operator in the present work is different. The condition $\delta h = -\frac{1}{2} d(\operatorname{tr}_g h)$ appeared naturally in the work of Chen and Nagano on harmonic metrics and harmonic symmetric 2-tensor fields, independently of its later applications in Ricci flow and Einstein geometry. Moreover, the formal adjoint operator $B_g^*$ naturally fits into the Berger–Ebin theory of the space of Riemannian metrics. Berger and Ebin [2] established fundamental injectivity results for a closely related operator differing only by the coefficient of the trace term. Analytically, the normalization adopted here agrees with the operator employed by Hamilton [7] in his study of the Ricci flow and with the coordinate-free Bianchi gauge formulation used by Kröncke [9].

The operator considered here has also appeared in earlier studies of harmonic diffeomorphisms and submersions, infinitesimal harmonic transformations, and Ricci solitons (see [15] and [16]). Thus, the present paper should be viewed as a continuation of a broader research program investigating the geometric role of the Chen-Nagano condition in global differential geometry. The novelty of the present work lies in employing this condition as a variational constraint for the Einstein-Hilbert functional.

**Remark 2.2.** The Chen-Nagano gauge condition $B_g(h) = 0$ formally coincides with the classical *de Donder harmonic gauge* (see [10]) used in general relativity and quantum gravity. However, the geometric interpretation is substantially different. In gravitational theory this condition is typically used as a gauge-fixing device for metric perturbations of Einstein's equations, whereas in the present paper it arises from the variational theory of Chen-Nagano harmonic tensors and serves as a constraint on admissible metric deformations. We emphasize that the de Donder gauge plays no role in the variational theory developed here and is mentioned only for comparison with the Chen-Nagano condition. Moreover, the Chen-Nagano condition may be interpreted as the infinitesimal harmonicity condition for the identity map $(M, g + h) \to (M, g)$. From this viewpoint, its appearance as the de Donder harmonic gauge reflects the underlying harmonic-map geometry of the Levi-Civita connection.

Define the vector space of Chen-Nagano harmonic tensors by

$$\mathcal{H}_g = \ker B_g = \left\{ h \in C^\infty(S^2M) : \delta h = -\frac{1}{2} d(\operatorname{tr}_g h) \right\}.$$

The operator $B_g$ considered here is conceptually equivalent to both the classical Killing operator and the conformal deformation operator studied by Berger and Ebin [2]. Since the injectivity of the principal symbol of $B_g^*$ is easily verified by a standard algebraic calculation in local coordinates, $B_g^*$ defines a first-order linear differential operator with injective principal symbol for $n \geq 2$ (see also [1, p. 456]; [2, p. 385]). This operator is called an *overdetermined elliptic operator*. Moreover, the associated second-order operator

$$B_g B_g^* = \frac{1}{2} \Delta_S + \frac{n-2}{4} d\delta$$

is strongly elliptic and formally self-adjoint, mirroring the foundational results established in [1, Appendix] and [2]. We recall that

$$\Delta_S\theta := 2\delta\, \delta^*\theta - d(\delta\theta)$$

is the Sampson Laplacian defined on 1-forms (see [17], [18]). In addition, $\ker B_g B_g^* = \ker B_g^*$ and $\dim\ker B_g B_g^* < \infty$.

Since $B_g^*\theta = \delta^*\theta + \frac{1}{2}(\delta\theta)g$ and $\mathrm{tr}_g(\delta^*\theta) = -\,\delta\theta,$ we obtain

$$| B_g^*\theta |^2 = | \delta^*\theta |^2 - (\delta\theta)^2 + \frac{n}{4}(\delta\theta)^2 = | \delta^*\theta |^2 + \frac{n-4}{4}(\delta\theta)^2.$$

Therefore, if $n \geq 4$ (resp. $n \leq 4$), then

$$| B_g^*\theta |^2 \geq | \delta^*\theta |^2 \quad \left(\text{resp. } | B_g^*\theta |^2 \leq | \delta^*\theta |^2\right).$$

In particular, $| B_g^*\theta |^2 = | \delta^*\theta |^2$ for $n = 4$.

More generally, for $n \geq 3$, taking the trace of $B_g^*\theta = \delta^*\theta + \frac{1}{2}(\delta\theta)g$ gives

$$\mathrm{tr}_g(B_g^*\theta) = \frac{n-2}{2}\delta\theta.$$

Hence $B_g^*\theta = 0$ implies $\delta\theta = 0$, and therefore $\delta^*\theta = 0$. Thus, for $n \geq 3$, the kernel of $B_g^*$ (resp. $B_g B_g^*$) coincides with the space of infinitesimal isometric transformations on $(M, g)$ (see [1, p. 40]).

On the other hand, by the classical Berger–Ebin–York decomposition theory based on Hodge–De Rham methods (see [1] and [2]), one obtains the $L^2$-orthogonal decomposition

$$C^\infty(S^2M) = \mathrm{Im}(B_g^*) \oplus_{L^2} \ker(B_g) = \mathrm{Im}(B_g^*) \oplus_{L^2} \mathcal{H}_g. \tag{2.1}$$

Consequently, every symmetric tensor $\varphi \in C^\infty(S^2M)$ admits a unique decomposition of the form

$$\varphi = B_g^*(\theta) + h,$$

where $h \in \mathcal{H}_g$ is a Chen-Nagano harmonic tensor.

Since the principal symbol of $B_g^*$ is injective but not surjective, the operator $B_g^*$ is overdetermined elliptic rather than classically elliptic. As a result, and in contrast to the elliptic case, the adjoint operator $B_g$ possesses an infinite-dimensional kernel (see Berger-Ebin [2]). Consequently, the space of Chen-Nagano harmonic tensors $\mathcal{H}_g = \ker B_g$ is infinite-dimensional (see also [3] and [4]).

Following Chen and Nagano, a symmetric tensor field $h \in C^\infty(S^2M)$ satisfying $B_g(h) = 0$ is *harmonic*. Correspondingly, variations generated by such tensors will be called *Chen-Nagano harmonic variations*. Therefore, the space $\mathcal{H}_g$ defines on $(M, g)$ a natural restricted class of Chen-Nagano symmetric harmonic variations of the metric $g$.

Thus, $\mathcal{H}_g = \ker B_g$ is the natural infinite-dimensional space of Chen-Nagano harmonic tensors and determines the admissible variational directions in the restricted variational problem considered below.

An additional decomposition will be useful for relating the Chen-Nagano condition to *infinitesimal harmonic transformations*.

Let

$$\mathcal{TT}_g(M) = \{h \in C^\infty(S^2M) : \delta h = 0,\ \mathrm{tr}_g\, h = 0\}$$

denote the infinite dimensional space of transverse-traceless tensors (see [1, Theorem 12.30]). By the classical Berger-Ebin-York decomposition theory, the algebraic sum

$$\operatorname{Im} \delta^* + C^\infty(M)g$$

is closed in $C^\infty(S^2M)$, and one has the $L^2$-orthogonal decomposition (see, for example, [1, pp. 118-119]; [8])

$$C^\infty(S^2M) = (\operatorname{Im} \delta^* + C^\infty(M)g) \oplus_{L^2} \mathcal{TT}_g(M).$$

Accordingly, every tensor $h \in C^\infty(S^2M)$ can be represented in the form

$$h = \delta^*\theta + \lambda g + h^{TT},$$

where

$$\theta \in C^\infty(T^*M), \qquad \lambda \in C^\infty(M), \qquad h^{TT} \in \mathcal{TT}_g(M).$$

Although the pair $(\theta,\ \lambda)$ need not be uniquely determined, the tensor

$$\delta^*\theta + \lambda g$$

is uniquely determined by the above orthogonal decomposition.

We recall that a smooth vector field $X$ on $(M, g)$ is called an *infinitesimal harmonic transformation* if the local one-parameter group of local diffeomorphisms generated by $X$ consists of harmonic diffeomorphisms (see [17]). Equivalently, if $\theta = X^\flat$ is the $g$-dual 1-form of $X$, then $X$ is an infinitesimal harmonic transformation if and only if $\Delta_S\theta = 0$.

**Proposition 2.1.** *Let $(M, g)$ be a compact Riemannian manifold of dimension $n \geq 3$. Suppose that*

$$h = \delta^*\theta + \lambda g + h^{TT},$$

*where $h^{TT} \in \mathcal{TT}_g(M)$. Then*

$$2B_g(h) = \Delta_S\theta + (n-2)\, d\lambda. \tag{2.2}$$

*Consequently, $h$ is a Chen-Nagano harmonic tensor if and only if*

$$\Delta_S\theta + (n-2)\, d\lambda = 0. \tag{2.3}$$

*In particular, any two of the following three assertions imply the remaining one:*

(i) *$h$ is a Chen-Nagano harmonic tensor;*

(ii) *$\theta^\#$ is an infinitesimal harmonic transformation;*

(iii) *$\lambda$ is constant.*

**Proof.** Since

$$\delta h^{TT} = 0, \qquad \mathrm{tr}_g\, h^{TT} = 0, \qquad \delta(\lambda g) = -\, d\lambda,$$

and

$$\mathrm{tr}_g(\delta^*\theta) = -\,\delta\theta,$$

we have

$$\delta h = \delta\delta^*\theta - d\lambda$$

and

$$d(\operatorname{tr}_g h) = -d(\delta\theta) + n\, d\lambda.$$

Therefore,

$$2B_g(h) \quad = 2\delta h + d\big(\operatorname{tr} h_g\big) = 2\delta\delta^*\theta - d(\delta\theta) + (n-2)d\lambda =$$

$$= \Delta_S\theta + (n-2)d\lambda.$$

This proves (2.2), while (2.3) follows from the definition

$$\mathcal{H}_g = \ker B_g.$$

The final assertion follows immediately from (2.2), since $n \geq 3$. □

**Corollary 2.2.** *Let $(M, g)$ be a compact Riemannian manifold with negative-definite Ricci tensor. Suppose that*

$$h = \delta^*\theta + \lambda g + h^{TT}, \quad h^{TT} \in \mathcal{TT}_g(M),$$

*is a Chen-Nagano harmonic tensor and that $\lambda$ is constant. Then*

$$\theta = 0.$$

*Consequently,*

$$h = \lambda g + h^{TT}.$$

*If, in addition,* h *is infinitesimally volume-preserving, namely*

$$\int_M tr_g\, h\, dv_g = 0,$$

*then $\lambda = 0$, and hence*

$$h = h^{TT}.$$

**Proof.** Since $h$ is Chen-Nagano harmonic and $\lambda$ is constant, equation (2.3) gives

$$\Delta_S\theta = 0.$$

Thus, $\theta^\sharp$ is an infinitesimal harmonic transformation. The Bochner-Yano formula for the Sampson Laplacian yields (see [17])

$$0 = \int_M \langle \Delta_S\theta, \theta\rangle\, dv_g = \int_M \big(|\nabla\theta^\sharp|^2 - \operatorname{Ric}_g(\theta^\sharp, \theta^\sharp)\big)\, dv_g.$$

Since $\mathrm{Ric}_g$is negative-definite, it follows that

$$\theta^\sharp = 0,$$

and therefore $\theta = 0$. Hence

$$h = \lambda g + h^{TT}.$$

Finally,

$$\mathrm{tr}_g\, h = -\delta\theta + n\lambda.$$

Integrating over $M$and using $\int_M \delta\theta \, d\, v_g = 0$, we obtain

$$\int_M \mathrm{tr}_g\, h\, dv_g = n\,\lambda\ \mathrm{Vol}(M, g).$$

Thus, the infinitesimal volume-preserving condition implies $\lambda = 0$, and consequently $h = h^{TT}$. □

## 3. Restricted Chen-Nagano Variational Principle

The purpose of this section is to develop a restricted variational theory for the Einstein-Hilbert functional based on Chen-Nagano variations. Unlike the classical Einstein variational principle, where admissible variations range over the entire Fréchet space $C^\infty(S^2M)$ of smooth symmetric tensor fields, we restrict the variational problem to the infinite-dimensional space $\mathcal{H}_g$ of Chen-Nagano harmonic tensors introduced in Section 2.

This restriction should be viewed as an infinite-dimensional analogue of a constrained extremum problem in classical variational calculus. The ambient Fréchet manifold of Riemannian metrics remains unchanged; only the admissible deformation directions are constrained by the Chen-Nagano gauge condition:

$$\delta\, h = -\frac{1}{2} d\big(\mathrm{tr}_g h\big),$$

for $h \in C^\infty(S^2M)$. As we shall see, this restriction nevertheless preserves a substantial part of the geometric content of the Einstein-Hilbert variational principle.

Let $S(g)$ denote the Einstein-Hilbert functional on the space of smooth Riemannian metrics on a compact Riemannian manifold $(M, g)$:

$$S(g) = \int_M s_g \ dv_g,$$

Let $g(t)$ be a smooth variation of $g$ and put $h = \frac{d}{dt}|_{t=0}\ g(t)$. The classical first variation formula (see Besse [1, §4.21]) is:

$$S_g'(h) = -\int_M \langle E_g, h\rangle dv_g, \tag{3.1}$$

where $E_g$ is the standard Einstein tensor:

$$E_g = Ric_g - \frac{1}{2} s_g g$$

To capture the constraint of volume preservation, we require that the smooth variation $g(t)$ does not alter the total volume of the manifold $M$. The corresponding infinitesimal deformation (or variational vector) $h$ must therefore satisfy the standard integral trace-free condition:

$$\int_M \mathrm{tr}_g \, h \, dv_g = 0. \tag{3.2}$$

Since the metric tensor $g$ itself belongs to the infinite-dimensional space $\mathcal{H}_g$ of Chen–Nagano harmonic tensors, condition (3.2) is equivalent to the global $L^2$-orthogonality relation: $h \perp g$. Consequently, the admissible volume-preserving Chen–Nagano variations are precisely characterized as symmetric tensor fields belonging to the intersection $h \in \mathcal{H}_g \cap g^{\perp}$.

To derive the corresponding Euler–Lagrange equation, consider the Berger–Ebin decomposition of the Einstein tensor:

$$E_g = B_g^*(\theta) + E_{g.}^H. \tag{3.3}$$

where $\theta$ is a smooth 1-form and $E_{g.}^H \in \mathcal{H}_g$. We shall refer to $E_{g.}^H$ as the *Chen-Nagano component* of the Einstein tensor.

The philosophy of the present construction differs essentially from other constrained variational approaches to the Einstein-Hilbert functional. In the theory of critical point metrics and related developments, one restricts the ambient space of metrics itself by imposing geometric conditions, such as constant scalar curvature. In contrast, the present approach preserves the full Fréchet manifold of Riemannian metrics and imposes a differential constraint only on the admissible deformation tensors. Thus, the restriction acts on tangent directions rather than on the underlying space of metrics.

The following theorem constitutes the main result of the paper.

**Theorem 3.1.** *Let $(M, g)$ be a compact Riemannian manifold of dimension $n \geq 3$. Then $g$ is a critical point of the Einstein-Hilbert functional with respect to all volume-preserving Chen-Nagano variations if and only if there exists a smooth 1-form $\theta$ such that:*

$$E_g = B_g^*(\theta) - \frac{n-2}{2n} \bar{s}\, g,$$

*where $\bar{s}$ is the average scalar curvature* of $(M, g)$.

**Proof.** Assume that $g$ is critical with respect to all admissible Chen-Nagano variations, meaning that the first variation (3.1) vanishes for all $\varphi \in \mathcal{H}_g \cap g^{\perp}$. Using decomposition (3.3), write:

$$E_g = B_g^*(\theta) + E_{g.}^H.$$

Substituting this decomposition into (3.1), we obtain:

$$\int_M \langle E_g, h\rangle \, dv_g = 0$$

for every $h \in \mathcal{H}_g \cap g^{\perp}$. Since the image Im $B_g^*$ is strictly orthogonal to $\mathcal{H}_g$ by the properties of the Berger–Ebin decomposition, it follows that:

$$\int_M \langle E_{g.}^H, h\rangle dv_g = 0$$

for every $h \in \mathcal{H}_g \cap g^{\perp}$. Hence, we have the geometric orthogonality:

$$E_{g.}^H \perp \left( \mathcal{H}_g \cap g^{\perp}\right).$$

Since the metric tensor $g$ itself belongs to $\mathcal{H}_g$, the orthogonal complement of $\mathcal{H}_g \cap g^{\perp}$ inside the subspace $\mathcal{H}_g$ is one-dimensional and generated precisely by $g$. Therefore, there must exist a constant $c \in \mathbb{R}$. such that: $E_{g.}^{H} = cg.$

Substituting this relation back into (3.3), we obtain the structural balance equation:

$$E_g = B_g^*(\theta) + c\,g \tag{3.4}$$

where $c = -\frac{n-2}{2n}\bar{s}$ for the *average scalar curvature* $\bar{s}$ of $(M, g)$ (see Remark 4.1).

Conversely, assume that equation (3.4) holds. Let $\varphi \in \mathcal{H}_g \cap g^{\perp}$ be an arbitrary admissible Chen-Nagano variation. By the orthogonality of the Berger–Ebin decomposition, we have:

$$\int_M \langle B_g^*(\theta), h\rangle \, dv_g = 0$$

while the volume-preserving condition yields:

$$\int_M \langle cg, h\rangle \, dv_g = 0$$

because $h \perp g$. Consequently, summing these equations gives:

$$\int_M \langle E_g, h\rangle \, dv_g = 0$$

and the classical variation formula (3.1) implies $\frac{d}{dt}|_{t=0}\, S(g(t)) = 0.$

Hence, $g$ is a genuine critical point of the Einstein-Hilbert functional with respect to all volume-preserving Chen-Nagano variations. The theorem is proved.

**Remark 3.1.** Theorem 3.1 shows that the Chen-Nagano gauge constraint does not destroy the variational structure of the Einstein-Hilbert functional. Although the admissible deformation directions are substantially restricted, the resulting Euler–Lagrange condition remains governed by the Einstein tensor. In this sense, the Chen-Nagano variational principle may be viewed as an infinite-dimensional analogue of a constrained extremum problem in classical variational calculus.

A noteworthy feature of this approach is that, unlike the theories of critical point metrics and related constrained variational frameworks, no restriction is imposed on the ambient space of Riemannian metrics itself. The entire Fréchet manifold of metrics is preserved, while the constraint acts only on tangent directions through the Chen-Nagano gauge condition. Remarkably, this weaker restriction is still sufficiently strong to yield the structural equation (3.4), which may be interpreted as a balance equation for the Einstein tensor.

**Definition 3.1.** *A Riemannian metric satisfying equation* (3.4) *will be called a Chen-Nagano critical metric*.

The terminology is justified by Theorem 3.1, which shows that equation (3.4) is precisely the Euler–Lagrange condition associated with the Einstein-Hilbert functional under volume-preserving Chen-Nagano variations. Thus, Chen-Nagano critical metrics play, within the restricted variational framework developed here, the same role that Einstein metrics play in the classical Einstein variational principle.

The significance of equation (3.4) extends beyond its variational origin. As will be shown in the next section, this structural equation naturally leads to generalized Ricci almost soliton-type geometries and, in the gradient case, to gradient Ricci almost solitons. Consequently, Chen-Nagano critical metrics provide a natural bridge between constrained variational geometry, orthogonal gauge decompositions, and modern Ricci flow theory.

## 4. Geometric Interpretation of Chen-Nagano Critical Metrics

The purpose of this section is to investigate the geometric and structural properties of Chen-Nagano critical metrics satisfying the constrained Euler–Lagrange equation derived in Theorem 3.1. We shall show that the structural equation obtained from the restricted variational principle naturally gives rise to geometric structures closely related to Ricci almost solitons and their gradient counterparts. In this way, the Chen-Nagano variational framework establishes a natural connection between constrained variational geometry and Ricci flow theory.

Recall from Definition 3.1 that a metric $g$ is called Chen-Nagano critical if there exist a smooth 1-form $\theta$ and a constant $c \in \mathbb{R}$ such that the Einstein tensor satisfies:

$$E_g = B_g^*\,(\theta)\ + c\,g \tag{4.1}$$

where $c = -\frac{n-2}{2n}\bar{s}$. By substituting the explicit definition of the Einstein tensor $E_g = Ric_g - \frac{1}{2}s_g g$ and the expression for the formal adjoint gauge operator $B_g^*(\theta) = \delta^*\theta + \frac{1}{2}(\delta\theta)g$ equation (4.1) can be explicitly rewritten as:

$$Ric_g - \frac{1}{2}s_g g = \delta^*\theta + \frac{1}{2}(\delta\theta)g + c\,g \tag{4.2}$$

Isolating the Ricci tensor on the left-hand side yields the following structural equation:

$$Ric_g = \delta^*\theta + \left(\frac{1}{2}\delta\theta +\ c + \frac{1}{2}s_g\right)g \tag{4.3}$$

Equation (4.3) may be regarded as the fundamental geometric form of the Chen-Nagano criticality condition. It shows that the Ricci tensor differs from the Einstein case by a gauge-generated deformation term and therefore provides a natural extension of the classical Einstein condition.

**Remark 4.1.** Taking the trace of equation (4.3), we obtain

$$s_g = -\,\delta\theta - \frac{2n}{n-2}c.$$

Integrating over the compact manifold $M$ and using the divergence theorem, we have

$$\int_M \delta\theta\, dv_g = 0.$$

Therefore, the total scalar curvature satisfies

$$S(g) = \int_M s_g\ dv_g = -\frac{2n}{n-2}c\,\mathrm{Vol}(M,g),$$

Consequently,

$$c = -\frac{n-2}{2n}\cdot\frac{S(g)}{\mathrm{Vol}(M,g)}$$

We recall that

$$\bar{s} = \frac{1}{\mathrm{Vol}(M)} \int_M s \, dv_g$$

is called the *average scalar curvature* of $(M, g)$ (see [1, pp. 114-117]). Then

$$c = -\frac{n-2}{2n} \bar{s}_g.$$

Thus, the constant c appearing in the structural equation is completely determined by the average scalar curvature of the metric and, in the unit-volume case, by the total scalar curvature. Therefore, equation (4.3) may be interpreted as a balance relation between the Einstein tensor, the Chen-Nagano gauge contribution $B_g^* (\theta)$, and the average (total, respectively) scalar curvature of $(M, g)$,

We now restrict our attention to the classical *gradient case*, where the generating 1-form $\theta$ is exact. Specifically, let $\theta = df$ for some smooth potential function $f \in C^\infty(M)$. In this setting, the symmetrized covariant derivative $\delta^*\theta$ reduces precisely to the Hessian operator of $f$, while the divergence $\delta\theta$ becomes the negative Laplacian of $f$: $\delta^*(df) = \mathrm{Hess}_g(f)$ and $\delta(df) = -\Delta_g f$.

Substituting these relations into the structural equation (4.3), we obtain:

$$Ric_g - \mathrm{Hess}_g(f) = \left(c + \tfrac{1}{2} s_g - \tfrac{1}{2} \Delta_g f\right) g \tag{4.4}$$

To relate equation (4.4) to the existing literature on geometric solitons, we introduce the smooth function $\lambda \in C^\infty(M)$ defined by:

$$\lambda = c + \tfrac{1}{2} s_g - \tfrac{1}{2} \Delta_g f. \tag{4.5}$$

for the Laplacian $\Delta_g f := \mathrm{tr}_g\left(\mathrm{Hess}_g f\right)$. Using this definition, equation (4.4) takes the elegant form:

$$Ric_g - \mathrm{Hess}_g(f) = \lambda \, g. \tag{4.6}$$

Equation (4.6) is precisely the defining relation for a *gradient Ricci almost soliton* (see, [11] and [12]). Unlike standard Ricci solitons where $\lambda$ must be a constant (see [7]), a Ricci *almost* soliton allows $\lambda$ to be a variable function on the manifold $M$. We have thus established the following fundamental geometric consequence:

**Theorem 4.1.** *Let* $(M, g)$ *be a Chen-Nagano critical metric. If the associated 1-form* $\theta$ *is exact,* $\theta = df$*, then the structural equation* (4.1) *determines a gradient Ricci almost soliton structure. The corresponding soliton function* $\lambda$ *is given explicitly by formula* (4.5).

Theorem 4.1 admits the following converse

**Theorem 4.2.** *Let* $(M, g)$ *be a compact Riemannian manifold of dimension* $n \geq 3$*, and suppose that* $(M, g, f, \lambda)$ *is a gradient Ricci almost soliton, that is,*

$$Ric_g - \mathrm{Hess}_g(f) = \lambda\, g.$$

*Then* $g$ *is a Chen–Nagano critical metric if and only if the soliton function* $\lambda$ *satisfies*

$$\lambda = \frac{1}{2}s_g - \frac{1}{2}\Delta_g f - \frac{n-2}{2n}\,\bar{s}_g,$$

*where*

$$\bar{s}_g = \frac{1}{Vol(M,g)}\int_M s_g\; dv_g$$

*is the average scalar curvature of* $(M, g)$.

**Proof.** Assume first that $g$ is a Chen–Nagano critical metric and that $\theta = df$. Then, by the structural equation,

$$E_g = B_g^*\,(df) - \frac{n-2}{2n}\bar{s}\; g.$$

Since

$$E_g = Ric_g - \frac{1}{2}\mathrm{s}_{\,g}\, g$$

and

$$B_g^*\,(df) = \mathrm{Hess}_g(f) - \tfrac{1}{2}\Delta_g f\; g$$

we obtain

$$Ric_g - \tfrac{1}{2}\mathrm{s}_{\,g}\, g = \mathrm{Hess}_g(f) - \tfrac{1}{2}\Delta_g f\; g - \tfrac{n-2}{2n}\bar{s}\; g.$$

Hence

$$Ric_g - \mathrm{Hess}_g(f) = \left(\tfrac{1}{2}\mathrm{s}_{\,g} - \tfrac{1}{2}\Delta_g f - \tfrac{n-2}{2n}\bar{s}\right) g. \;.$$

Therefore, the corresponding gradient Ricci almost soliton function is

$$\lambda = \frac{1}{2}s_g - \frac{1}{2}\Delta_g f - \frac{n-2}{2n}\bar{s}.$$

Conversely, suppose that $(M, g, f, \lambda)$ is a gradient Ricci almost soliton and that

$$\lambda = \frac{1}{2}s_g - \frac{1}{2}\Delta_g f - \frac{n-2}{2n}\bar{s}.$$

Then

$$Ric_g - \mathrm{Hess}_g(f) = \left(\frac{1}{2}s_g - \frac{1}{2}\Delta_g f - \frac{n-2}{2n}\bar{s}\right) g.$$

Rearranging this equality gives

$$Ric_g - \frac{1}{2}s_g\, g = \mathrm{Hess}_g(f) - \frac{1}{2}\Delta_g f\, g - \frac{n-2}{2n}\bar{s}\, g.$$

Since $B_g^*(df) = \mathrm{Hess}_g(f) - \frac{1}{2}\Delta_g f\, g$ we obtain

$$E_g = B_g^*(df) - \frac{n-2}{2n}\bar{s}\, g.$$

By Theorem 3.1, $g$ is a critical point of the Einstein–Hilbert functional with respect to all volume-preserving Chen–Nagano variations. Thus $g$ is a Chen–Nagano critical metric. The proof is complete. □

On the other hand, it is worth noting that the explicit functional form of the variable soliton function $\lambda = c + \frac{1}{2}s_g - \frac{1}{2}\Delta_g f$ establishes a direct conceptual link with the structural trace identities investigated in [13] for infinitesimal transformations under geometric flows. This agreement underscores the fact that the Chen-Nagano restricted variational principle naturally recovers the analytical environment of soliton rigidity within a self-contained optimization framework.

Theorem 4.1 shows that Chen-Nagano critical metrics are naturally related to geometric structures arising in Ricci flow theory (see [6] and [7]). In particular, the restricted variational principle developed in Section 3 does not merely produce a formal Euler–Lagrange equation. Rather, it leads to a broader class of geometric structures governed by equations closely connected with gradient Ricci almost solitons (see [11] and [12]).

Thus, the Chen-Nagano variational principle extends the classical Einstein variational framework in a nontrivial way. While Einstein metrics arise when the gauge term vanishes, the general Chen-Nagano criticality condition allows a wider family of geometric structures. The next section shows, however, that under additional geometric assumptions this apparent flexibility becomes strongly restricted and the classical Einstein condition re-emerges through a global rigidity phenomenon.

## 5. Global Rigidity of Chen-Nagano Critical Metrics

In this section, let $(M, g)$ be a compact Chen-Nagano critical metric of dimension $n \geq 3$. Thus, there exists a smooth 1-form $\theta$ such that

$$E_g = B_g^*(\theta) + cg, \tag{3.4}$$

where

$$c = -\frac{n-2}{2n}\,\bar{s}$$

is the average scalar curvature of $(M, g)$.

The first result shows that the gauge 1-form occurring in the restricted Euler-Lagrange equation belongs to the kernel of the Sampson Laplacian pointwise.

**Proposition 5.1**. *Let* $(M, g)$ *be a compact Chen-Nagano critical metric satisfying* (3.4). *Then*

$$\Delta_S\,\theta = 0. \tag{5.1}$$

*Equivalently, the vector field*

$$X = \theta^\sharp$$

*is an infinitesimal harmonic transformation.*

**Proof**. Write equation (3.4) in the form

$$E_g = \delta^*\theta + \mu g,$$

where

$$\mu = \frac{1}{2}\delta\theta + c.$$

Applying Proposition 2.1 to this decomposition, we obtain

$$2B_g(E_g) = \Delta_S\theta + (n-2)d\mu. \tag{5.2}$$

On the other hand, the contracted second Bianchi identity gives

$$\delta E_g = 0,$$

whereas

$$\mathrm{tr}_g\, E_g = -\frac{n-2}{2} s_g.$$

Therefore,

$$2B_g(E_g) = 2\delta E_g + d(\mathrm{tr}_g\, E_g) = -\frac{n-2}{2} ds_g. \tag{5.3}$$

Taking the trace of (3.4), we find

$$-\frac{n-2}{2} s_g = \frac{n-2}{2}\delta\theta + n\, c.$$

Hence,

$$s_g = -\delta\theta - \frac{2n}{n-2} c,$$

and consequently

$$d\mu = \frac{1}{2} d(\delta\theta) = -\frac{1}{2} ds_g. \tag{5.4}$$

Substituting (5.4) into (5.2) and comparing with (5.3), we get

$$-\frac{n-2}{2} ds_g = \Delta_S\theta - \frac{n-2}{2} ds_g.$$

Therefore,

$$\Delta_S\theta = 0.$$

The proof is complete. □

**Corollary 5.1.** *Under the assumptions of Proposition* 5.1*, the following energy identity holds:*

$$\int_M |\, \delta^*\theta\, |^2\ dv_g = \frac{1}{2}\int_M (\delta\theta)^2\, dv_g. \tag{5.5}$$

**Proof.** Since

$$\Delta_S\theta = 0,$$

integration by parts gives

$$0 \quad = \int_M \langle\, \Delta_S\theta, \theta\rangle\, dv_g = 2\int_M |\, \delta^*\theta\, |^2\ dv_g - \int_M (\delta\theta)^2\, dv_g$$

.

This is precisely (5.5). □

The preceding proposition gives an immediate rigidity theorem under a sign condition on the Ricci curvature.

**Corollary 5.2**. *Let* $(M, g)$ *be a compact Chen-Nagano critical metric of quasi-negative Ricci curvature. Then*

$$\theta = 0,$$

*and* $g$ *is an Einstein metric.*

**Proof.** Here *quasi-negative Ricci curvature means* that $\mathrm{Ric}_g \leq 0$ on $M$, and that for some point $x \in M$ the bilinear form $\mathrm{Ric}_{g,x}$ is negative definite.

Let $X = \theta^\sharp$. The standard Bochner-Yano integral formula for the Sampson Laplacian gives (see [17])

$$\int_M \langle \Delta_S \theta, \theta \rangle \, dv_g = \int_M \big( | \nabla X |^2 - \mathrm{Ric}_g(X, X) \big) \, dv_g. \tag{5.6}$$

By Proposition 5.1, the left-hand side of (5.6) vanishes. Since the Ricci curvature is non-positive, both terms in the right-hand side are non-negative. Hence

$$\nabla X = 0.$$

Moreover, quasi-negative Ricci curvature means that $\mathrm{Ric}_g$is negative definite at least at one point. Therefore $X$vanishes at that point and, because $X$ is parallel, it follows that

$$X = 0$$

on $M$. Consequently,

$$\theta = 0.$$

Equation (3.4) now reduces to

$$E_g = cg.$$

Taking the trace, we obtain

$$-\frac{n-2}{2} s_g = nc,$$

and hence

$$\mathrm{Ric}_g = \frac{s_g}{n} g.$$

Therefore $g$ is Einstein. □

We next obtain a stronger rigidity theorem under the classical assumption of constant scalar curvature. In contrast with the previous version of the result, no exactness assumption on the gauge 1-form is required.

**Theorem 5.1**. *Let $(M, g)$ be a compact Chen-Nagano critical metric of dimension $n \geq 3$. If the scalar curvature $s_g$ is constant, then $g$ is Einstein. Moreover,*

$$B_g^*(\theta) = 0, \qquad \delta\theta = 0, \qquad \delta^*\theta = 0.$$

*Consequently, the vector field $\theta^\sharp$ is Killing. If $\theta = df$, then $f$ is constant.*

**Proof.** Since $s_g$ is constant,

$$B_g(E_g) = \delta E_g + \frac{1}{2} d(\mathrm{tr}_g\, E_g) = -\, \frac{n-2}{4}\, ds_g = 0.$$

Applying $B_g$ to equation (3.4) and using

$$B_g(cg) = 0,$$

we obtain

$$B_g B_g^*(\theta) = 0.$$

Taking the $L^2$-inner product with $\theta$, we find

$$0 \quad = \int_M \langle B_g B_g^*(\theta), \theta\rangle \, dv_g = \int_M |\, B_g^*(\theta)\, |^2 \; dv_g.$$

Hence,

$$B_g^*(\theta) = 0.$$

It follows from (3.4) that

$$E_g = cg,$$

and thus $g$ is Einstein.

Taking the trace of

$$B_g^*(\theta) = 0,$$

we obtain

$$0 \quad = \mathrm{tr}_g\, B_g^*(\theta) = \mathrm{tr}_g \left(\delta^*\theta + \frac{1}{2}(\delta\theta) g\right) = \frac{n-2}{2} \delta\theta$$

Therefore,

$$\delta\theta = 0.$$

Consequently,

$$\delta^*\theta = 0.$$

Since

$$\delta^*\theta = \frac{1}{2}\mathcal{L}_{\theta^\sharp} g,$$

the vector field $\theta^\sharp$ is Killing.

Finally, if

$$\theta = df,$$

then

$$\delta(df) = -\Delta_g f = 0.$$

Since $M$ is compact, $f$ is constant. □

The next consequence remains useful when the scalar curvature is not assumed to be constant.

**Corollary 5.3**. *Let $(M, g)$ be a compact gradient Chen-Nagano critical metric of dimension $n \geq 3$. Assume that*

$$\theta = df,$$

*and let*

$$X = \nabla f.$$

*If*

$$X(s_g) \geq 0$$

*on $M$, then $f$ is constant and $g$ is Einstein.*

**Proof.** Applying $B_g$ to equation (3.4), using

$$B_g(cg) = 0,$$

and recalling that

$$B_g(E_g) = -\frac{n-2}{4} ds_g,$$

we obtain

$$B_g B_g^*(df) = -\frac{n-2}{4} ds_g. \tag{5.7}$$

Taking the $L^2$-inner product with $df$, we get

$$\int_M |B_g^*(df)|^2 \, dv_g = -\frac{n-2}{4}\int_M \langle d\, s_g, df\rangle \, dv_g = -\frac{n-2}{4}\int_M X(s_g)\, dv_g.$$

The left-hand side is non-negative, whereas the right-hand side is non-positive by assumption. Therefore, both sides vanish, and hence

$$B_g^*(df) = 0.$$

Taking the trace gives

$$0 = \mathrm{tr}_g\, B_g^*(df) = \frac{n-2}{2}\delta(df).$$

Thus,

$$\Delta_g f = 0.$$

By compactness, $f$ is constant. Equation (3.4) then reduces to

$$E_g = cg,$$

which proves that $g$ is Einstein. □

**Remark 5.2.** Proposition 5.1 shows that the Chen-Nagano gauge 1-form arising from the restricted Einstein-Hilbert variational principle is automatically an infinitesimal harmonic transformation. Thus, the constrained Euler-Lagrange equation imposes a substantially stronger restriction on the gauge term than is visible from the earlier integral identity alone. Under natural curvature assumptions, the apparent flexibility of the Chen-Nagano criticality equation collapses to the classical Einstein condition.

## 6. Negative Einstein Metrics and Restricted Stability

In this section, we examine the restricted second variation of the Einstein-Hilbert functional at a compact Einstein metric with negative Einstein constant (see (1)). We shall prove that, in this case, the Chen-Nagano gauge condition together with infinitesimal preservation of the scalar curvature automatically reduces to the classical transverse-traceless gauge. Consequently, this allows us to study the variational stability of such metrics, which is defined by the nonpositivity of the second variation of the functional on the space of admissible variations (see [1]; [2]; [8]).

Let $(M, g)$ be a compact Einstein manifold satisfying

$$\mathrm{Ric}_g = \lambda g, \qquad \lambda < 0.$$

Recall that a symmetric 2-tensor $h \in C^\infty(S^2T^*M)$ is a Chen-Nagano harmonic tensor if

$$\delta h = -\frac{1}{2} d(\mathrm{tr}_g\, h).$$

We denote by

$$\mathcal{H}_g^{sc} = \{h \in \mathcal{H}_g\colon\ s_g'(h) = 0\}$$

the space of infinitesimally scalar-curvature-preserving Chen-Nagano variations.

We first show that, in the negative Einstein case, these two spaces coincide.

**Proposition 6.1.** *Let $(M, g)$ be a compact Einstein manifold satisfying*

$$\mathrm{Ric_g} = \lambda \mathrm{g}, \lambda < 0.$$

Then

$$\mathcal{H}_\mathrm{g}^{\mathrm{sc}} = \mathcal{TT}_\mathrm{g}(\mathrm{M}).$$

*Equivalently, every infinitesimally scalar-curvature-preserving Chen-Nagano variation is transverse-traceless.*

**Proof.** Let $h \in \mathcal{H}_g^{\mathrm{sc}}$, and put

$$\tau := \mathrm{tr}_g\, h.$$

In the sign convention of Besse, the first variation of the scalar curvature is given by

$$s_g'(h) = \Delta\tau + \delta(\delta h) - \langle \mathrm{Ric}_g, h\rangle. \tag{6.1}$$

Since $h$satisfies the Chen-Nagano gauge condition,

$$\delta h = -\frac{1}{2} d\tau,$$

we obtain

$$\delta(\delta h) = -\frac{1}{2}\delta d\tau = -\frac{1}{2}\Delta\tau.$$

Moreover, since $g$ is Einstein,

$$\langle \mathrm{Ric}_g, h\rangle = \lambda\tau.$$

Therefore, equation (6.1) becomes

$$s_g'(h) = \frac{1}{2}\Delta\tau - \lambda\tau. \tag{6.2}$$

By assumption, $s_g'(h) = 0$. Hence

$$\Delta\tau = 2\lambda\tau. \tag{6.3}$$

Multiplying equation (6.2) by $\tau$, integrating over $M$, and using Green's formula, we obtain

$$0 = \frac{1}{2}\int_M |\, d\tau\,|^2 \; dv_g - \lambda \int_M \tau^2 \; dv_g. \tag{6.4}$$

Since $\lambda < 0$, both terms in the right-hand side of (6.4) are nonnegative, and the second one is strictly positive unless $\tau \equiv 0$. Therefore,

$$\mathrm{tr}_g\, h = \tau = 0.$$

The Chen-Nagano condition now yields

$$\delta h = -\frac{1}{2} d\tau = 0.$$

Thus $h \in \mathcal{TT}_g(M)$.

Conversely, if $h \in \mathcal{TT}_g(M)$, then the Chen-Nagano condition is immediate. Moreover,

$$s_g'(h) = -\langle \mathrm{Ric}_g, h\rangle = -\lambda\, \mathrm{tr}_g\, h = 0.$$

Hence $h \in \mathcal{H}_g^{\mathrm{sc}}$. This proves $\mathcal{H}_g^{\mathrm{sc}} = \mathcal{TT}_g(M)$. □

We now apply this result to the restricted second variation. Recall that, at an Einstein metric, the second variation formula may be written as

$$S_g''(h,h) = -\frac{1}{2}\int_M \langle\, \Delta_E h, h\rangle \, dv_g$$

for $h \in \mathcal{H}_g^{sc} = \mathcal{TT}_g(M)$, where $\Delta_E$ is the *Einstein elliptic operator* on symmetric 2-tensors (see Besse [1, p. 132] and Koiso [8]).

This gives the following stability criterion.

**Corollary 6.1.** *Let* $(M, g)$ *be a compact Einstein manifold with*

$$Ric_g = \lambda g, \qquad \lambda < 0.$$

*Assume that*

$$\int_M \langle \Delta_E h, h\rangle \, dv_g \geq 0$$

*for every* $h \in \mathcal{TT}_g(M)$. *Then*

$$S_g''(h,h) \leq 0$$

*for every infinitesimally scalar-curvature-preserving Chen--Nagano variation* $h$. *In particular,* $g$ *is variationally stable with respect to scalar-curvature-preserving Chen-Nagano variations if and only if the Einstein elliptic operator is nonnegative on* $\mathcal{TT}_g(M)$.

**Remark 6.1.** The negativity of the Einstein constant does not by itself determine the sign of the quadratic form in (6.6). Its role is different and more precise: it eliminates all trace directions compatible with the Chen-Nagano gauge and the infinitesimal scalar-curvature-preserving condition. Consequently, the restricted variational problem reduces exactly to the classical spectral stability problem for the Lichnerowicz operator on transverse-traceless tensors.

For comparison, when $\lambda = 0$, equation (6.3) implies that $\operatorname{tr}_g h$ is constant; this remaining trace component disappears after imposing the infinitesimal volume-preserving condition. When $\lambda > 0$, nontrivial trace directions may occur precisely when $2\lambda$ belongs to the spectrum of the Laplacian on functions.

**Remark 6.2.** The nonnegativity assumption on the Einstein operator in Corollary 6.1 is satisfied for several standard classes of negative Einstein manifolds. In particular, if $(M, g)$, $n \geq 3$, is a compact hyperbolic manifold, then the Einstein operator $\Delta_E$ is positive on $\mathcal{TT}_g(M)$. Moreover, if $(M, g, J)$ is a compact Kähler--Einstein manifold with nonpositive Einstein constant, then $\Delta_E$ is nonnegative on $\mathcal{TT}_g(M)$. In the negative Kähler-Einstein case, its kernel may be nontrivial and is related to infinitesimal deformations of the complex structure (see [5]).

**Data Availability Statement**

No datasets were generated or analyzed during the current study.

**Funding Declaration**

The author declares that no financial funds, grants or other support were used in the preparation of this manuscript.

**Conflict of interest statement**

The author declares no conflict of interest.

Sergey Stepanov

Department of Mathematics All Russian Institute for Scientific and Technical Information of the Russia Academy of Sciences
20, Usievicha street, Moscow, Russia 125190

Department of Mathematics and Data Analysis, Finance University
49-55, Leningradsky Prospect, Moscow, Russia 125468
e-mail: s.e.stepanov@mail.ru

Irina Tsyganok

Department of Mathematics and Data Analysis, Finance University
49-55, Leningradsky Prospect Moscow Russia 125468
e-mail: i.i.tsyganok@mail.ru